\documentclass[12pt]{amsart}
\usepackage{graphicx}
\usepackage{amscd}
\usepackage{amsmath}
\usepackage{amsfonts}
\usepackage{amssymb}
\usepackage{verbatim}

\newcommand{\toplines}{\hline \hline \\[-1.06ex]}
\newcommand{\sepline}{\\[-2.13ex] \hline \\[-2.13ex]}

\newcommand{\botlines}{\\[+0.8ex] \hline \hline}

\newcommand{\prob}[1]{\ensuremath{{\rm P}\!\left( #1 \right)}}

\newcommand{\ber}[1]{\textit{Ber\ensuremath{\mspace{2mu}(#1)}}}
\newcommand{\firule}[1]{\varphi_{\mathtt{#1}}}
\newcommand{\jfi}[1]{\mathtt{#1}}
\newcommand{\jf}{\mathtt{j}}
\newcommand{\Supp}[1]{\mathrm{Supp}{(#1)}}

\newtheorem{theorem}{Theorem}
\theoremstyle{plain}

\newtheorem{proposition}{Proposition}

\newcommand{\x}{x}
\newcommand{\X}{X}
\newcommand{\PN}[1]{{\mathbb P\!}_N\!\left({#1}\right)}
\newcommand{\PI}[1]{{\mathbb P\!}_*\!\left({#1}\right)}
\newcommand{\PIkaal}{{\mathbb P\!}_*}
\newcommand{\prcI}[2]{\ensuremath{{\mathbb P\!}_*\!}\left( #1 \, |\, #2\right)}
\newcommand{\prcN}[2]{\ensuremath{{\mathbb P\!}_N\!}\left( #1 \, |\, #2\right)}
\newcommand{\vecphi}{\varphi}
\newcommand{\map}{\Psi_{\varphi}}
\newcommand{\Map}{\Psi_{\!\Phi}}
\newcommand{\e}[1]{\varepsilon_{#1}}

\begin{document}
\DeclareGraphicsExtensions{.jpg,.pdf,.png}

\parindent0pt

\title[Cellular Automata]
{Stability in random Boolean cellular automata on the integer lattice}

\author{F. Michel Dekking, Leonard van Driel and Anne Fey}
\address{ Delft Institute of Applied Mathematics,
 Technical University of Delft,  Mekelweg 4, 2628 CD Delft, The Netherlands\\
 email: F.M.Dekking@tudelft.nl}
\address{  Eurandom, Eindhoven, The Netherlands \\
 email: Fey@eurandom.tue.nl}

\maketitle
\begin{abstract}
We consider one-dimensional random boolean cellular automata,
i.e., the cells are identified with the integers from 1 to  $N$.
The behavior of the automaton is mainly determined by the support
of the random variable that selects one of the sixteen possible
Boolean rules, independently for each cell. A cell is said to stabilize if
 it will not change its state anymore  after some time.
 We classify the one-dimensional random boolean automata according to the positivity of their probability of
 stabilization. Here is an example of a consequence of our results: if the
support contains at least 5 rules, then asymptotically as
$N\rightarrow \infty$ the probability of stabilization is positive, whereas
 there exist random boolean cellular automata
 with 4 rules in their support for which this probability tends to  0.
\end{abstract}

\section{Introduction}

We study the dynamics of cellular automata as amply discussed in
Stephen Wolfram's book ``A New Kind of Science" (\cite{NKS}).
The cells, which are identified with the integers from 1 to $N$, can be in two
possible states 0 or 1, and the state at the next time instant is a function of
the states of the two neighbouring cells. For this  function, also called (evolution) rule,
there are 16 different possibilities, see also Table \ref{tab-bca-phi}.
Although there is a chapter in \cite{NKS} on \emph{random} initial conditions for
CA's, there is hardly anything on \emph{random} rules, or even on the possibility to
assign different rules to different cells (so called \emph{inhomogeneous} CA's).
In this paper we will give a classification of one-dimensional Random Boolean CA's.

Our results are also of interest in the context of Random Boolean Networks
 introduced by Stuart Kauffman in 1969 (\cite{Kauffman1}). Here not only the rules
 and the initial condition are chosen randomly, but also the neighbours of  a cell,
 i.e., the $N$ cells obtain a random graph structure.
For each cell $i$ independently a mother $M(i)$ and a father $F(i)$ are chosen
uniformly from the remaining pairs of cells. This gives a graph
with vertices $i$ and directed edges $(M(i),i)$, $(F(i),i)$ for
$i=1,\dots,N$. Next, each cell $i$ is assigned independently a
Boolean rule $\Phi_i$ uniformly from the set of 16 Boolean rules
(i.e., maps from $\{0,1\}^2$ to $\{0,1\}$, see also Table
\ref{tab-bca-phi}). To start, all cells obtain independently an
initial value $X_i(0)$ from the set $\{0,1\}$ with equal
probability. The state of a RBN  at time $t$ will be
denoted by $\X(t)$. Then the dynamics is given by the map
$\Map$:
$$\X(t+1)=\Map(\X(t))$$
 defined  by
$$X_i(t+1) = \Phi_i(X_{M(i)}(t),X_{F(i)}(t))$$
for all cells $i$.

\begin{figure} [h]
\begin{center}
\includegraphics[width=14cm,clip= bbllx=0, bblly=0, bburx=300, bbury=160]{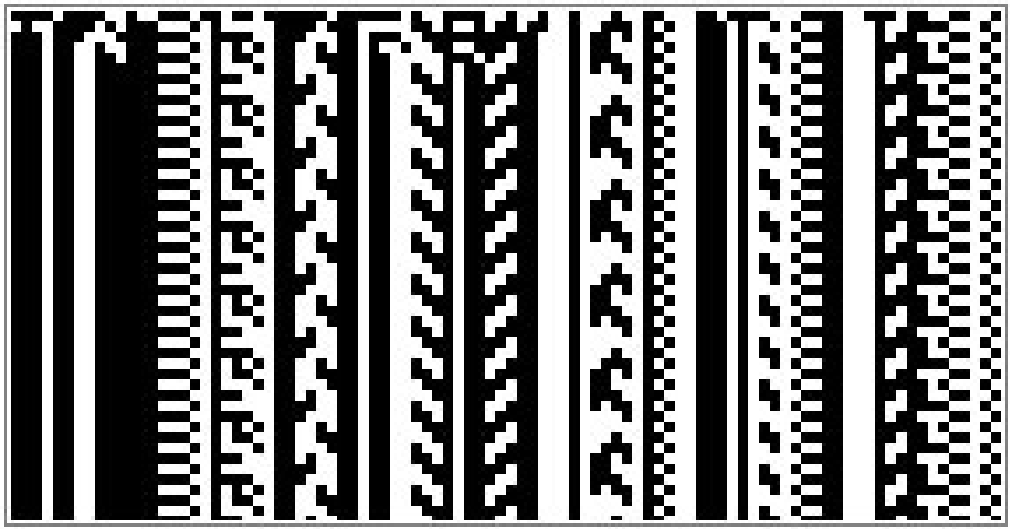}\
\caption{\footnotesize Evolution of a one-dimensional RBCA.
  After two time steps, most cells already
become periodic. After five steps, the whole CA is periodic.}
\end{center}\label{fig-intro}
\end{figure}

A $1$-dimensional random Boolean cellular automaton can be considered as an RBN with $$M(i)\equiv i-1\quad\mathrm{ and } \quad F(i)\equiv i+1,$$
for $i\!=\!2,\dots,N\!-\!1$. We put $M(1)\!=\!N, F(1)\!=\!2$, and $M(N)\!=\!N-1, F(N)\!=\!1$.
In Figure \ref{fig-intro} we display a realization of the evolution of a RBCA in the
usual way: zero's and ones are coded by white and black squares,
space is in the horizontal direction, time in the vertical
direction.

We say a cell \emph{stabilizes} if its  value is not changed anymore
after some time. In Figure~\ref{fig-intro} the stable cells yield
the black or white stripes; in this RBCA example about 56\% of the
cells stabilize.

Kauffman (\cite{MR762660}) found in his
simulations that about 60-80\% of the cells in a two-connected RBN stabilize. It was
later shown rigorously by Luczak and Cohen (\cite{luczak:1991})
that in the limit as the number of cells $N$ tends to infinity
100\% of the cells will stabilize---size $N=10\,000$ in Kauffmans
simulations is much too small to unearth this behavior.\\
For the RBN's there is no natural model for the
infinite size ($N\rightarrow\infty$) network. For a one-dimensional CA
 there is: the graph is simply the set of all
integers with edges between nearest neighbors. It therefore makes
sense to speak
 of the probability $\sigma_*$ of stabilization of the infinite
 automaton. The result of Luczak and Cohen leads to the question
 whether this $\sigma_*=1$. In this paper we will show that this
 is rarely the case, and moreover that there are many cases where
 actually  $\sigma_*=0$.
 By different cases we mean different distributions on the set of 16 Boolean
 rules. We remark that the results in \cite{MR762660} and
\cite{luczak:1991} are restricted to the uniform distribution on
the set of rules---see \cite{Drossel-2007,Lynch-1993,Lynch-1993-2,Lynch-PhD} for
some results for RBN's with other rule distributions.\\
The second problem we will discuss is whether
$\sigma_N\rightarrow\sigma_*$ when $N\rightarrow\infty$. We will
show that this will often be the case, but may also fail to hold.\\
 For some previous work with rather partial
results on 2--dimensional versions of  RBCA's see
e.g.,  \cite{fogelman}. We will give an application of our results
at the end of  Section~\ref{imp}  to a recent paper on RBN's of Greil and Drossel (\cite{Drossel-2007}).

Finally we mention that our interest in RBCA was also motivated by
certain CA models  of the subsurface in the geosciences
(see ~\cite{geo}).

\section{Random Boolean cellular automata}\label{HD}

Before we define RBCA, we
describe their realizations which are inhomogeneous cellular
automata.
 The dynamics in such a  CA is defined by attaching a Boolean rule to each cell.
The rule is given by specifying what its value is on the four
different values
 $(1,1)$, $(1,0)$, $(0,1)$ and $(0,0)$ of the neighbors of a cell.
 In Table \ref{tab-bca-phi}, all 16 different rules
are defined.
\begin{table}[h!]\centering
\caption{ Definition of the 16 Boolean rules.  The last column
gives a formula for  $\firule{j}(x,y)$ when $\firule{j}$ is
affine.}
\begin{tabular}{|l|c c c c|c|}
\hline $\jf$ & $\!\firule{j}(0,0)\!$ &\!$\firule{j}(0,1)$\!
    &\! $\firule{j}(1,0)$\! &\!$\firule{j}(1,1)$ & affine\phantom{$I^{I^2}_{I_2}$} \\
\hline
0 &0 & 0 & 0  &0 &   $0$    \\
1  &1 & 0 & 0  &0 &   no     \\
2  &0 & 1 & 0  &0 &   no     \\
3  &1 & 1 & 0  &0 &   $x+1$  \\
4  &0 & 0 & 1  &0 &   no     \\
5  &1 & 0 & 1  &0 &   $y+1$  \\
6  &0 & 1 & 1  &0 &   $x+y$  \\
7  &1 & 1 & 1  &0 &   no     \\
8  &0 & 0 & 0  &1 &   no     \\
9  &1 & 0 & 0  &1 &   $x+y+1$\\
10 &0 & 1 & 0  &1 &   $y$    \\
11 &1 & 1 & 0  &1 &   no     \\
12 &0 & 0 & 1  &1 &   $ x$   \\
13 &1 & 0 & 1  &1 &    no    \\
14 &0 & 1 & 1  &1 &    no    \\
15 &1 & 1 & 1  &1 &    $1$   \\
\hline
\end{tabular}
\label{tab-bca-phi}
\end{table}

For ease of notation these are indexed as $\firule{j}$ where $\jf$
is the number in binary notation obtained from the four bits in
the row describing  $\firule{j}$. Anytime we refer to one of these
16 rules we use index $\jf$ (note the different font); when we
use index $i$ we refer abusively to the rule assigned to cell $i$.

An inhomogeneous CA's mapping $\map: \{0,1\}^N\to \{0,1\}^N$ is
determined by a vector of rules,  $\vecphi=(\varphi_1, \ldots,
\varphi_N)$. Here for all time steps, $\varphi_i$ is the rule of
cell $i$. Hence given a state $\x = (x_1 \ldots x_N)$ we get
$(\map(\x))_i = \varphi_i(x_{i-1},x_{i+1})$ for $i=2,\dots, N-1$.
To deal with the problem that occurs in the end points we consider
the automata as defined on a circle; a convenient way to express
this is to define $x_{i+kN}=x_i$ for $i=1,\dots, N$, and all
integers $k$.

To address  $\map$  we will also write
the vector of lower indices of the $\varphi_i$ $(1 \leq i \leq N)$
mentioned. E.g., mapping (0,3,3,15) describes
$\map:\{0,1\}^4\to\{0,1\}^4$ with rule vector
$(\firule{0},\firule{3},\firule{3},\firule{15})$.\\

In a RBCA the rules $\varphi_i$ are chosen according to a
distribution $\Phi$ on the set of rules, independently for each
$i$. The initial state is an i.i.d. sequence $(X_i(0))$ from a
Bernoulli(1/2) distribution, independent from the $(\Phi_i)$. Then
the value of the $i^\textrm{th}$ cell at time $t+1$ is given by
\begin{equation}
X_i(t+1)=\Phi_i(X_{i-1}(t),X_{i+1}(t)).
\end{equation}
We denote by ${{\mathbb P\!}_N}$ the probability measure generated by the  ($X_i(0)$) and
the  ($\Phi_i$) on the space of infinite periodic sequences with period $N$.
It is natural to consider also the probability measure ${{\mathbb P\!}_*}$ on $\{0,1\}^{\mathbb{Z}}$
generated by two infinite independent i.i.d. sequences ($X_i(0)$) and ($\Phi_i$).

\section{Stability}\label{stab}

We say
cell $i$ \emph{stabilizes at time} $T$ if
$$X_i(t)=X_i(T) \quad \text{for all } t\ge T,$$
and $T$ is the smallest integer with this property. Let $\sigma_N$
be the probability that a cell in a random cellular automaton of
size $N$ with rule distribution $\Phi$ stabilizes. Because of
stationarity,
$$\sigma_N=\PN{\text{cell } 0 \text{ stabilizes}}.$$

We first give an example were $\sigma_N$ does not converge as
$N\rightarrow\infty$.\\
Consider the RBCA's for $N=1,2\dots$ with rule distribution given by
$$\PN{\Phi=\firule{6}}=1.$$
Note that randomness is only involved because of the randomness in
the initial state $\X(0)=(X_i(0))$, and that $\firule{6}$ is the
well known linear rule which can be written as
$\firule{6}(x,y)=x+y$. This addition, as all additions  with 0's
and 1's, is modulo 2.

\begin{proposition}\label{prop:discontin}
For the  $\firule{6}$-RBCA one has
$$ \liminf_{N\rightarrow\infty} \sigma_N=0,
    \qquad \limsup_{N\rightarrow\infty} \sigma_N =1.$$
\end{proposition}

\textit{Proof:} It is well known (see e.g.~\cite{algebra}) and
easy to prove that the $\firule{6}$-automaton satisfies
$$X_i(t)=X_{i-t}(0)+X_{i+t}(0) \qquad \text{ if }  t=2^p \text{ for
some } p.$$
It follows that $\sigma_N=1$ for each $N=2^p$ which is a power of 2:
$$X_i(N)=X_{-N+i}(0)+X_{N+i}(0)=2X_i(0)=0 \qquad \text{for all } i,$$
which implies that $X_i(t)=0$ for all  $t\ge N$. This proves that
$\limsup \sigma_N =1$. To obtain the statement on $\liminf$,
consider the subsequence of $N$ of the form $N=2^p+1$. Then we find for all $i$
$$X_i(N-1)=X_i(2^p)=X_{-2^p+i}(0)+X_{2^p+i}(0)=X_{i+1}(0)+X_{i-1}(0). $$
From this we directly obtain that
$$X_0(N)=X_{-1}(N-1)+X_{1}(N-1)=X_{-2}(0)+X_{0}(0)+X_{0}(0)+X_{2}(0),$$
and so
\begin{equation}\label{even}
X_0(N)=X_{-2}(0)+X_{2}(0).
\end{equation}
Thus
$$\PN{X_0(N-1)\ne X_0(N)}=\PN{X_{-1}(0)+X_{1}(0)\ne X_{-2}(0)+X_{2}(0)}=\frac12.$$
Note that this computation does not only hold for $t=N-1$, but for any $t$ which
satisfies $t=-1 \mod N$.
An infinite sequence of such $t$'s is given by $t=2^{(2k+1)p}, k=0,1,2,\dots$.
For $N$ of the form  $2^p+1$ we therefore obtain that
$$1-\sigma_N=\PN{X_0(t)\ne X_0(t+1)\quad \text{for infinitely many  } t}\ge\frac12.$$
This inequality is a warm up for a more general result: we can  generalize
Equation~(\ref{even}) to  for  all $m$
$$X_0(N+2^m-2)=X_{-2^m}(0)+X_{2^m}(0).$$
Since for $m=0,1,\dots,\ell$ the $X_{-2^m}(0)+X_{2^m}(0)$ are independent Bernoulli
 variables with success probability $1/2$, it follows that
$$\PN{\exists m, 0\le m < \ell: X_0(N+2^m-2)\ne X_0(N-1)}=1-\frac1{2^\ell},$$
and employing the same infinite sequence of $t$'s as for the $\ell=1$ case,
we can deduce that  $ \liminf_{N\rightarrow\infty} \sigma_N=0$.\quad$\Box$ \\

It is convenient to define  random variables $Z_i$ by $Z_i=1$ if
cell $i$ stabilizes, and 0 otherwise. We then have (independently of
$i$):
$$\sigma_N=\PN{Z_i=1}, \quad \text{and} \quad  \sigma_*=\PI{Z_i=1}.$$
The \emph{uniform} RBCA is given by the rule distribution
$\PN{\Phi=\firule{j}}=1/16$ for $j=0,1,\ldots,15$.

\begin{proposition} \label{unif}
Let  $\sigma_N$ be the probability of stabilization of a size $N$ uniform RBCA,
$\sigma_*$ that of an infinite size uniform RBCA. Then
$$\lim_{N\rightarrow\infty}\sigma_N=\sigma_*\quad \textrm{and}\quad \sigma_*>0.$$
\end{proposition}

\textit{Proof:} The two rules $\firule{0}$ and $\firule{15}$ are
called \emph{walls}. We denote $W=\{\firule{0}$,$\firule{15}\}$.
Cells that have obtained a rule from $W$ are already stable at
time 0 or 1. We define for positive integers $k$ and $l$ the
\emph{chambers} $C_{k,l}$ by
 $$C_{k,l}=\{\Phi_{-k}\in W, \Phi_{-k+1}\notin W, \Phi_{-k+2}\notin W,\ldots,\Phi_{\l-1}\notin W,\Phi_l\in W\}.$$
In addition we put $C_{0,0}=\{\Phi_0\in W\}$. Since the
chambers are disjoint events whose union has probability 1 we can
write
 \begin{eqnarray*}
 \sigma_*=\PI{Z_0=1}=\prcI{Z_0=1}{C_{0,0}}\PI{C_{0,0}}+
        \sum_{k,l=1}^{\infty}\prcI{Z_0=1}{C_{k,l}}\PI{C_{k,l}}.
\end{eqnarray*}
We will make a similar splitting for finite RBCA of size $N$. In case
$N=2M+1$ we consider the indices of the chambers
 from the  set  $\{-M,\dots,0,\dots,M\}$, in case
$N=2M$ from $\{-M\!+\!1,\dots,0,\dots,M\}$. In the sequel we will
assume $N$ is odd, as our arguments can be transferred trivially to
the case $N$ even. Let $D$ be the event
$$D=\Big(C_{0,0}\cup\bigcup_{k,l=1}^M C_{k,l}\Big)^c.$$
Then $D$ occurs iff there is no wall among $\Phi_{-M}, \dots, \Phi_0$ or among
$\Phi_{0}, \dots, \Phi_M$.
Since $\PN{\Phi_i\notin W}=\PI{\Phi_i\notin W}=7/8$, we obtain that
$$\PI{D}=\PN{D}=2\Big(\frac78\Big)^{M+1}\!\!-\Big(\frac78\Big)^{2M+1}
    \le2\Big(\frac78\Big)^{N/2}.$$
The crux of the proof is that for all $1\le k,l\le M$ and for
$k=l=0$ we have
$$\prcN{Z_0=1}{C_{k,l}}=\prcI{Z_0=1}{C_{k,l}},$$
simply because both the finite automaton and the infinite
automaton evolve between the walls independently of the evolution outside the walls.
Thus
\begin{eqnarray*}
 |\sigma_N-\sigma_*|&\le & \prcI{Z_0=1}{D}\PI{D}+\prcN{Z_0=1}{D}\PN{D}\\
          &\le & 4\Big(\frac78\Big)^{N/2}.
\end{eqnarray*}
It follows that $\sigma_N \rightarrow \sigma_*$  exponentially fast.
Note that
$$\sigma_N\ge\PN{\Phi\in W}=\frac18,$$
and so $\sigma_*\ge 1/8>0$.  \hfill $\Box$\\

 Actually we can obtain a good estimate of $\sigma_*$ using the exponential
 convergence of $\sigma_N$, by sampling uniform RBCA's for some large $N$ and computing a
 95\% confidence interval. We found that $\sigma_{100}=0.678058\pm 0.000029$.
 Interestingly, an exhaustive enumeration for $N=6$ yields already a value $\sigma_6$
 close to 0.68.

 We end this section with the remark that the `wall'-property directly implies that
 the sequence $(Z_i)$ is strongly mixing under $\PIkaal$:
 $$\PI{Z_0=\e{0}, Z_i=\e{1}}\rightarrow \PI{Z_0=\e{0}}\PI{Z_0=\e{1}}\quad \mathrm{as} \quad i\rightarrow\infty.$$
  Therefore the ergodic theorem applies, yielding that we
 also have pathwise about 68\% of stable cells.

\section{Stability by means of impermeable blocks}\label{imp}
Note that  the proof of Proposition~\ref{unif} will still be valid if
$\Phi$ is not uniformly distributed, but has an arbitrary
distribution with  $\PN{\Phi=\firule{j}}>0$  for all $\jf$ (of
course $7/8$ has to be replaced by
$\PN{\Phi=\firule{0}}+\PN{\Phi=\firule{15}}$). The only thing that
matters is the \emph{support} of $\Phi$ defined by
$$\text{Supp}(\Phi)=\{\jf: \PN{\Phi=\firule{j}}>0\}=\{\jf: \PI{\Phi=\firule{j}}>0\}.$$
In this section we will generalize the concept of a wall. An easy
adaptation of its proof  will permit to draw the  conclusions of
Proposition~\ref{unif} for a large majority of the $2^{16}$ different
supports.

Consider an inhomogeneous BCA. We call adjacent cells  $\{i+1,
\ldots,i+p\}$ an \emph{impermeable block} of size $p$, if there
exist  $b=(b_1,\ldots,b_p)\in\{0,1\}^p$, and
$(\jfi{j_1},\ldots,\jfi{\jfi{j_p}})\in\{0,\ldots,15\}^p$ such that
if
$$(x_{i+1}(0),\ldots,x_{i+p}(0))=(b_1,\ldots,b_p)\quad \text{and}
\quad
(\varphi_{i+1},\ldots,\varphi_{i+p})=(\varphi_{\jfi{j_1}},\ldots\varphi_{\jfi{j_p}}),$$
then for all $t$
$$ (x_{i+1}(t),\ldots,x_{i+p}(t)), $$
is unchanged, whatever values are chosen for $x_i(s)$ and
$x_{i+p+1}(s)$, for $s=0,1\dots,t-1$.
Here we assume that $N$ is larger than $p+2$.\\
As an example take $p=2$, $b=(0,0)$ and  $(\jfi{j_1},\jfi{j_2})=(2,4)$.\\
If a particular $\Phi$ has  $\firule{j}$'s in its support that
yield an impermeable block, then any $\Phi$ having this set of $\firule{j}$'s  in
its support will also have an impermeable block. We therefore look
for minimal sets of $\jf$'s such that the $\firule{j}$'s give rise
to impermeable blocks. Here is a listing of these sets:
\begin{eqnarray*}
\mathcal{G}=\big\{\{0\}, \{1\}, \{7\}, \{8\}, \{14\}, \{15\}, \{2,
4\}, \{2, 5\},
\{2, 9\}, \{2, 12\}, \{2, 13\},\quad\\
\{3, 4\},\{3, 5\}, \{3,10\},  \{3, 13\}, \{4, 9\}, \{4, 10\}, \{4, 11\}, \{5,11\},\quad\\
\{5,12\},\{6, 11\},\{6, 13\}, \{10, 12\}, \{10, 13\}, \{11,12\}, \{11, 13\},\quad\\
 \{2, 6, 10\}, \{4, 6, 12\}, \{9, 10, 11\},\{9, 12, 13\}\big\}.
\end{eqnarray*}
Because of the symmetries of the collection of sets of CA rules
(cf.~Section~\ref{symm}), of these 30 subsets there is a
subcollection $\mathcal{\widetilde{G}}$ of only 12 where the
dynamics of the associated $\Map$ is essentially different.
Table~\ref{imperm} gives these supports  with their
$\varphi$-blocks
$(\varphi_{\jfi{j_1}},\ldots\varphi_{\jfi{j_p}})$, and their
$b$-blocks $(b_1,\ldots,b_p)$.

\begin{table}[h!]\centering
\caption{Impermeable blocks.}
 {
 \begin{tabular}{lcc||lcc}
 \toplines
 Support &  $\varphi$-block  & $b$-block & \quad Support &  $\varphi$-block  & $b$-block\\
 \sepline
 $\{0\}$ \quad     & (0)   & (0)          &\quad $\{2,12\}$\quad  & (2,12)   & (0,0)  \\
 $\{1\}$ \quad     & (1,1,1)   & (0,1,0)  &\quad $\{2, 13\}$\quad  &  (13, 2)  &  (1,0) \\
 $\{8\}$ \quad     & (8,8)     & (0,0)    &\quad $\{3, 5\}$ \quad  & (5,3)     &  (0,0)   \\
 $\{2,4\}$ \quad  &  (2,4)    & (0,0)    &\quad $\{3,10\}$ \quad  & (10,3)    & (0,0)    \\
 $\{2,5\}$ \quad  &  (5,2)    & (1,0)    &\quad $\{10,12\}$\quad  & (10,12)   &  (0,0)   \\
 $\{2,9\}$\quad  &(2,9,9,2)  &(0,0,1,0) &\quad $\{2,6,10\}$\quad & (10,6,2)  & (1,1,0)
\botlines\label{imperm}
\end{tabular}}
\end{table}

Often there is more than one  $b$-block that will work, in the table we
have chosen  $b$-blocks with minimal length. It can be quickly
verified that the entries in Table~\ref{imperm} yield impermeable
blocks. See Appendix~\ref{imp-ver} for two examples of the type of
verification one has to perform.

We can describe the supports which do \emph{not} give rise to impermeable blocks by the list
\begin{eqnarray*}
\mathcal{B}=\big\{ \{2, 3, 6\}, \{2, 3, 11\}, \{2, 10, 11\},\{3,
9, 11\}, \{4, 5, 6\}, \{4, 5, 13\}, \quad\\  \hspace*{5cm} \{4,
12, 13\}, \{5, 9, 13\}, \{3,6, 9, 12\},
 \{5,6, 9, 10\} \big\}.
\end{eqnarray*}
Any support which does not admit an impermeable block is a subset
of a set in $\mathcal{B}$. For an example of a proof that with
these supports no impermeable blocks occur, see
Appendix~\ref{no-imp}.

One can check that any subset of $\{0,1,\ldots,15\}$ is either a
subset of a set from $\mathcal{B}$, or contains a set from
$\mathcal{G}$, and so we have completely described the supports
with regard to the existence of impermeable blocks.

\medskip

{\bf Example}. In \cite{Drossel-2007} the two-connected RBN with solely rule $\firule{7}$ is studied,
and it is claimed that ``every node oscillates with period 2". Since $\firule{7}$ is mirror equivalent to $\firule{1}$, we see that actually the RBCA with support only consisting of $\firule{7}$, has
absorbing block (1,0,1). Let $i$ be a node in the $\firule{7}$-RBN  with parents $j$ and $k$. Then, if $i$ happens to be one of the parents of both $j$ and $k$ (figure 8 configuration), and the state of $i$ is 0, and that of $j$ and $k$ is 1, then these three nodes keep that state forever. This gives that the expected number of stable nodes in the RBN is at least
$$N\cdot \PN{\mathrm{figure\; 8\; configuration}}\cdot\PN{\mathrm{initial\; values\; 1,0,1}}=N\left(\frac2{N-1}\right)^2.\frac18.$$

\section{Classifying RBCA's}\label{class}

In the previous section we have established that all the RBCA's whose support
contains a set in $\mathcal{G}$ are well behaved in the sense that they are stability
continuous ($\sigma_N\rightarrow\sigma_*$), and do have a positive
probability to stabilize ($\sigma_*>0$). As there are only 53 supports described
by $\mathcal{B}$ this means that  this regular behaviour  manifests itself for at
least 99,9\% of the RBCA's. The following result gives a complete classification.

\begin{theorem}\label{main}
An infinite one-dimensional random Boolean cellular automaton has no stable cells, i.e., $\sigma_*=0$, if and only if the support $\Supp{\Phi}$ of the rule distribution $\Phi$ is a subset of $\{\firule{2}, \firule{10}\}, \{\firule{10}, \firule{11}\}, \{\firule{4}, \firule{12}\}, \{\firule{12}, \firule{13}\},
\{ \firule{3}, \firule{6}, \firule{9}, \firule{12}\}$ or $\{ \firule{5}, \firule{6}, \firule{9}, \firule{10}\}$.
\end{theorem}

\textit{Proof:} From the results of the previous section it is clear that any  CA with $\sigma_*=0$ must have its support in (the subsets of ) collection $\mathcal{B}$.
Because of symmetries it is enough to consider the 18 supports which occur as subsets of
the sets in the collection
$\mathcal{\widetilde{B}}$ given by
\begin{eqnarray*}
\mathcal{\widetilde{B}}=\big\{\{2, 3, 6\}, \{2, 3, 11\}, \{2, 10,
11\}, \{3,6, 9, 12\}\big\}.
\end{eqnarray*}
In Section~\ref{linbadguys} we will prove that  subsets of  rules in
$\{3,6,9, 12\}$ all give supports that belong to CA with $\sigma_*=0$.

It occurs that the remaining supports can be split into two parts.
The first part consists of the supports which occur as  subsets of
$\mathcal{\widetilde{B}}$ given by
\begin{eqnarray*}
\mathcal{\widetilde{S}}=\big\{\{2, 3\}, \{2,6\},\{2, 11\}, \{2, 3,6\},\{2,
3,11\},\{2,10,11\}\big\}.
\end{eqnarray*}
The $\widetilde{}$\, here indicates that we removed those supports
that are symmetry equivalent to supports that are already listed,
as, for example, $\{3,11\}$. The supports in
$\mathcal{\widetilde{S}}$ will be shown in Section~\ref{absorb}
 to be of the same regular type as  found in Section~\ref{imp}.
The second part consists of the two supports
$$\{2\} \textrm{ and} \; \{2,10\}.$$
Here the automata behave very similarly.
\begin{figure}[b]
\begin{center}
\includegraphics[width=13cm,clip= bbllx=0, bblly=0, bburx=403, bbury=202]{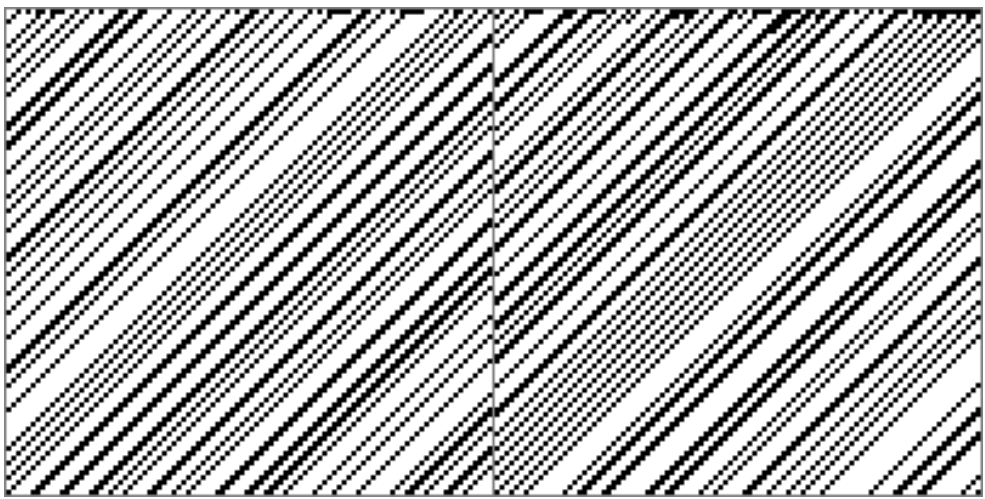}
\caption{\footnotesize Evolution of the RBCA with support $\{\varphi_2\}$ and $\{\varphi_2,\varphi_{10}\}$ .}
\end{center}\label{fig-2-10}
\end{figure}
First note that $\firule{10}$ is simply the left shift. The crucial fact is that
$$ \firule{2}(x,y)=\firule{10}(x,y) \quad\mathrm{if}\quad (x,y)\ne (1,1).$$
First we consider the case $\Supp{\Phi}=\{\firule{2}\}$. Then $\firule{2}$ acts as
the left shift for all times $t\ge 1$, since after one application of $\firule{2}$
there will be no blocks $1  \ast 1$ (i.e., 101 or 111) left anymore, since in the
initial sequence $0\ast 1$ cannot overlap with itself with lag 2. The map  $\firule{2}$
acts then as the shift on a subshift in which blocks of more than two ones can not occur.
It follows that for each $i$ the process $(X_i(t)_{t\ge 1})$ is the same subshift,
which clearly can not stabilize. (Actually it can be shown that this process is a
hidden Markov chain on 8 equiprobable states.)

In case $\Supp{\Phi}=\{\firule{2}, \firule{10}\}$, the CA runs as a left shift with occasional transformations of 1's into 0's. This is a form of self-organization. In fact any block of 1's that occurs at time 0 to the right of a cell $i$ with $\Phi_i=\firule{2}$ will be cut into blocks of 1 and 11 separated by 0's.
Since any cell $i$ has a cell $i'$ to the right of it with $\Phi_{i'}=\firule{2}$ with probability 1,
there exists a (random) $t_0$, such that the process $(X_i(t))_{t\ge t_0}$ will only contain
 1 and 11, and cannot stabilize. \hfill $\Box$

 \section{Stability by means of absorbing blocks}\label{absorb}

Here we show that we may define a weaker notion than that of an
impermeable block which still permits us to conclude to continuity
of stabilization and a positive probability to stabilize.\\
Consider an inhomogeneous CA. We call adjacent cells  $\{i+1,
\ldots,i+p\}$ an \emph{absorbing block} of size $p$, if there
exist  $b$ in $\{0,1\}^p$, and   $(\jfi{j_1},\ldots,\jfi{j_p})$
in $\{0,\ldots,15\}^p$ such that  if
$$(x_{i+1}(0),\ldots,x_{i+p}(0))=(b_1,\ldots,b_p)\quad \text{and}
\quad
(\varphi_{i+1},\ldots,\varphi_{i+p})=(\varphi_{\jfi{j_1}},\ldots\varphi_{\jfi{j_p}}),$$
then for some $c$ with $1<c<p$ the value $x_{i+c}(t)$ is unchanged
for all $t$, whatever values are chosen for $x_i(s)$ and
$x_{i+p+1}(s),s=0,1\dots,t-1$.
(Here we assume that $N$ is larger than $p+2$).
The idea is that the evolution of one or both adjacent cells $i$ and $i+p+1$
may `penetrate' the block, but will never influence the `central'
cell $i+c$.\\
We give an example for an inhomogeneous automaton with rules
$\firule{2}$ and $\firule{11}$. Let $p=4$, $b=(0,1,1,0)$, and
$(\jfi{j_1},\ldots,\jfi{j_4})=(2,2,11,2)$. Then this block is
absorbing with $c=2$:
\begin{center}
  \begin{tabular}{c|c|c|c|c|c}
cell $i$  & cell $i+1$  &  cell $i+2$&  cell $i+3$ &  cell $i+4$ &  cell $i+5$\\
\;  &  $\firule{2}$ &   $\firule{2}$ & $\firule{11}$ & $\firule{2}$ & \; \\
 \sepline
 $x$   & 0      & 1 & 1 & 0 & $y$ \\
$x'$   & $1-x$  & 1 & 0 & 0 & $y'$\\
$x''$  & $1-x'$ & 0 & 0 & $y'$ & $y''$\\
$x'''$ & 0      & 0 & 1 & $y''$ & $y'''$\\
       & 0      & 1 & 1 & 0 &
\end{tabular}
\end{center}
Note that at time $t=4$ the block $b=(0,1,1,0)$ reappears, so cell
$i+2$ has a periodic evolution with period 4, independently of the
evolution of cell $i$ and cell $i+5$. This establishes that
$(2,2,11,2)$ is an absorbing block. Similarly it can be verified
that $(2,2,3)$ is an absorbing block for the automata with rules
$\firule{2}$ and $\firule{3}$.

For the support $\{2,6\}$, we present an absorbing block \emph{family}
with stable cells: take $\varphi$-block
$$(2,2,2,6,6,6,2,2,2,6,6,6,2,2,2,2)$$ and $b$-block family
$$(0,0,1,1,0,1,0,1,0,1,1,0,x,y,z,w),$$ where $x,y,z$ and $w$ can
have any value as long as both $xz$ and $yw$ are not $11$. As a
consequence of the $\firule{2}$ rule, neither $xz$ nor $yw$ can be
$11$ at any later time. One can then check, in a case by case
analysis for all remaining possibilities for $xz$ and $yw$ (three
possibilities each), that the $b$-block family keeps reappearing,
and that the fourth and seventh cell are stable. We find that if
$yw = 01$ or $yw = 10$, it takes 8 time steps for another member
of the $b$-block family to reappear, for all other possibilities,
this takes 4 time steps.

It follows from these observations that  the rules with supports in the collection
$\mathcal{\widetilde{S}}$ considered in the proof of Theorem~\ref{main} do not have stable cells.

We can not yet conclude that $\sigma_*>0$ for supports $\{\firule{2},\firule{3}\}$ and $\{\firule{2},\firule{11}\}$ , since in both cases the
central cell has a period 4 dynamics. However, we can find larger
absorbing blocks with a stable central cell $c$: for $\{2,3\}$
take $p=9$, $c=6$, $\varphi$-block: $(2,2,3,2,3,2,2,2,3)$ and
$b$-block: $(0,0,1,1,0,0,0,0,1)$. For $\{2,11\}$ take $p=10$,
$c=6$, $\varphi$-block: $(2,2,11,2,11,2,2,2,11,2)$ and  $b$-block:
$(0,1,1,0,0,0,0,1,1,0)$. The existence of such absorbing blocks
with a stable cell implies that $\sigma_*>0$.


\section{Affine chaos}\label{linbadguys}

Here we will analyze the RBCA's whose support is a subset of $\{3,6,9,12\}$.

\begin{figure} [h]
\begin{center}
\includegraphics[width=13cm,clip= bbllx=0, bblly=0, bburx=402, bbury=202]{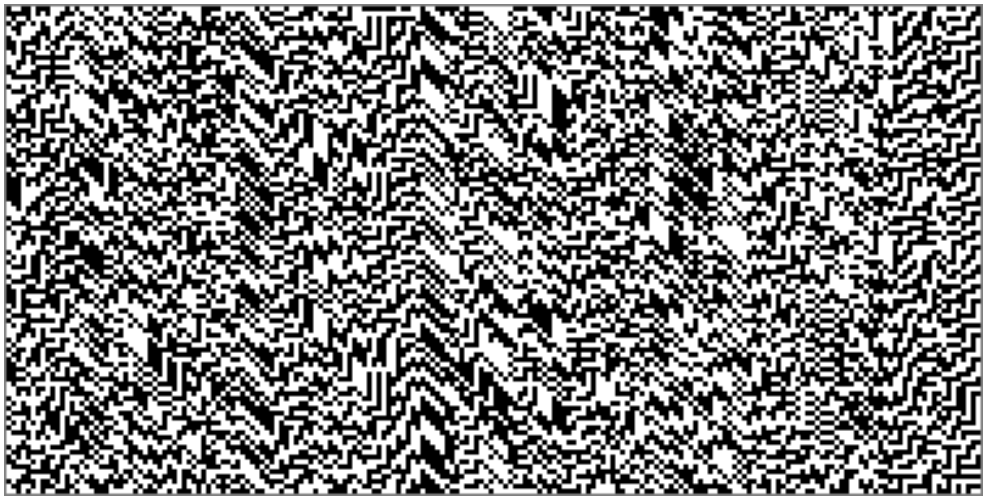}
\caption{\footnotesize Evolution of an RBCA with support $\{3,6,9,12\}$.}
\end{center}\label{fig-36912}
\end{figure}

The crucial observation  is that not only these rules are affine,
but they can all be written as
$$\firule{j}(x,y)=x+f(y,\firule{j}),$$
where $f(y,\firule{j})=1,\,y,\,y+\!1,\,0$ for $\jf=3,6,9,12$ respectively. It
follows from this  that we can write for some function $f_t$ with $4t-2$ arguments
\begin{equation}\label{xprop}
x_i(t)=x_{i-t}(0)+f_t(x_{i-t+2}(0),\dots,x_{i+t}(0),\varphi_{i-t+1},\dots,\varphi_{i+t-1}).
\end{equation}
Indeed, for $t=1$ we have
$$x_i(1)=x_{i-1}(0)+f_1(x_{i+1}(0),\varphi_i),$$
where $f_1(y,\varphi_i)=f(y,\varphi_i )$ as above.
For $t=2$ we obtain, applying this equation or its shift three times
\begin{eqnarray*}
x_i(2)&=&x_{i-1}(1)+f_1(x_{i+1}(1),\varphi_i)\\
      &=&x_{i-2}(0)+f_1(x_{i}(0),\varphi_{i-1})+
            f_1(x_{i}(0)+f_1(x_{i+2}(0),\varphi_{i+1}),\varphi_i).
\end{eqnarray*}
Defining
\begin{eqnarray*}
 f_2(x_{i}(0),\dots,x_{i+2}(0),\varphi_{i-1},\dots,\varphi_{i+1})\hspace*{5cm}\\
\hfill = f_1(x_{i}(0),\varphi_{i-1})+f_1(x_{i}(0)+f_1(x_{i+2}(0),\varphi_{i+1}),\varphi_i),
\end{eqnarray*}
we obtain Equation~(\ref{xprop}) for $t=2$.
Continuing in this fashion we obtain  Equation~(\ref{xprop}) for
all $t$.

 The following proposition will immediately imply that $\sigma_*=0$ for
 the RBCA's from this section.

\begin{proposition}
Consider an RBCA with $\mathrm{Supp}(\Phi) \subset
\{\firule{3},\firule{6},\firule{9},\firule{12}\}$, and let $i$ be
a fixed cell. Then the random variables
$(X_i(t))_{t=0}^{t=\infty}$ are independent \ber{\frac12} random
variables under $\PIkaal$.
\end{proposition}

\textit{Proof:}
We prove this by induction w.r.t.~the length $t$ of the cylinders. We will show that
for each cell $i$ and for all $t=1,2,\dots$
\begin{equation}\label{ber}
\PI{X_{i}(0)=\e{0},X_i(1)=\e{1},\dots,X_i(t-1)=\e{t-1}}=2^{-t}
\end{equation}
for all $(\e{0},\e{1},\dots,\e{t-1})$ from$\{0,1\}^{t}$.

\begin{figure} [h]
\begin{center}
\includegraphics[height=8cm]{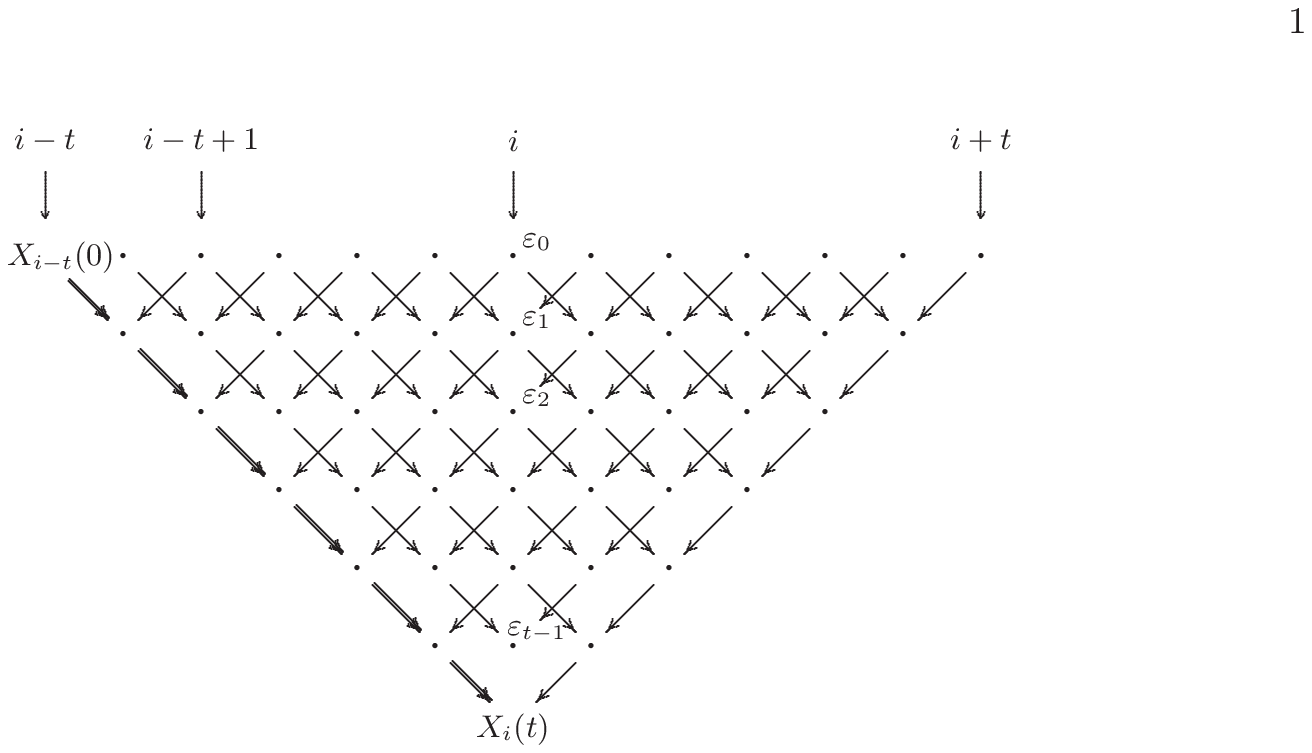}
\caption{\label{fig-prop2} \footnotesize Additive propagation.}
\end{center}
\end{figure}


Although in general the $(X_i(t))$ will not be stationary (any  cell which stabilizes
but not at $t=0$ has a transient evolution!), a simple computation shows that Equation (\ref{ber})
implies that shifted cylinders of length $t$ also have probability $2^{-t}$ to occur.

We will use the following abbreviations:
\begin{eqnarray*}
     Y_{i,t}       &=& f_t(X_{i-t+2}(0),\dots,X_{i+t}(0),\Phi_{i-t+1},\dots,\Phi_{i+t-1})\\
E_{i,t,\varepsilon} &=& \{X_{i}(0)=\e{0},X_i(1)=\e{1},\dots,X_i(t-1)=\e{t-1}\}.
\end{eqnarray*}
First note that Equation~(\ref{ber}) is true for $t=1$ by
definition of the RBCA. Figure~\ref{fig-prop2} is useful in the
remainder of the proof.


 Suppose that Equation~(\ref{ber}), i.e.,
$\PI{E_{i,t,\varepsilon}}=2^{-t}$ has been
 proved for cylinders of length $t$.
    Then for all $\e{0},\dots,\e{t-1}$ from$\{0,1\}^{t}$
\begin{eqnarray*}
&&\PI{E_{i,t,\varepsilon},X_i(t)=1}=\\
&&\PI{E_{i,t,\varepsilon},X_{i-t}(0)=0,Y_{i,t}=1}+
 \PI{E_{i,t,\varepsilon},X_{i-t}(0)=1,Y_{i,t}=0}=\\
&&\PI{X_{i-t}(0)=0}\PI{E_{i,t,\varepsilon},Y_{i,t}=1}+
 \PI{X_{i-t}(0)=1}\PI{E_{i,t,\varepsilon},Y_{i,t}=0}=\\
&&\tfrac12\big(\PI{E_{i,t,\varepsilon},Y_{i,t}=1}+
 \PI{E_{i,t,\varepsilon},Y_{i,t}=0}\big)=\\
&&\tfrac12\PI{E_{i,t,\varepsilon}}=2^{-(t+1)},
\end{eqnarray*}
where we used that $X_{i-t}(0)$ is independent of
$E_{i,t,\varepsilon}$ and of $Y_{i,t}$. This equality proves
Equation~(\ref{ber}) for length
$t+1$ cylinders.\hfill \qed\\

What remains for the  automata in this section is the question
whether $\sigma_N \rightarrow \sigma_{*}$. For some automata this
is easy to answer. We have already seen that for support $\{6\}$
there is \emph{no} stability continuity (see Proposition
\ref{prop:discontin}). By symmetry, the same is true for support
$\{9\}$ and with a little more work this can also be shown for
support $\{6, 9\}$. On the other hand it is easy to show that
$\sigma_{2N+1}=0$ and $\sigma_{2N}=2^{-N}$  for $\mathrm{Supp}(\Phi)=\{3\}$ and
that $\sigma_N=2^{-N+1}$ for $\mathrm{Supp}(\Phi)=\{12\}$, hence we do
have stability continuity in these cases. Also for the combined
support $\{3, 12\}$ it is easily shown that $\sigma_N\le
2^{-N+1}$, hence here too there is stability continuity. However,
for the automata with supports $\{3,6\}$, $\{6,12\}$, $\{3,6,9\}$,
$\{3,6,12\}$, $\{6,9,12\}$, and $\{3,6,9,12\}$, the computation of
$\sigma_N$ becomes quite involved.  We conjecture that
there is stability continuity for these cases.

\section{Conclusion}

We have studied the behavior of one-dimensional random Boolean cellular automata with two
inputs. Although this behavior can be quite
diverse, we have shown that it does not depend on the
probabilities with which the Boolean rules
$\varphi_{\jfi{0}},\dots,\varphi_{\jfi{15}}$ are attached to the
integers, but only on the support of this random variable $\Phi$,
i.e. on the probabilities $\prob{\Phi=\varphi_{\jf}}$ being 0 or
positive. This contrasts with the behavior of random Boolean networks as
determined by Lynch (\cite{Lynch-PhD, Lynch-2006}),
 see also \cite{Drossel-2007}. His result is
that there is ordered behavior as long as
$$\prob{\Phi=\varphi_{\jfi{0}}}+\prob{\Phi=\varphi_{\jfi{15}}}\ge
\prob{\Phi=\varphi_{\jfi{6}}}+\prob{\Phi=\varphi_{\jfi{9}}},$$ and
chaotic behavior when the opposite inequality holds. Nevertheless,
there is an interesting parallel: for the one-dimensional random cellular automata
$\prob{\Phi=\varphi_{\jfi{0}}}+\prob{\Phi=\varphi_{\jfi{15}}}>0$
implies regular behavior in the sense that $\sigma^*>0$ and
exponential convergence of $\sigma_N$ to $\sigma^*$, while
$\prob{\Phi=\varphi_{\jfi{6}}}+\prob{\Phi=\varphi_{\jfi{9}}}=1$
implies that $\sigma^*=0$, and \emph{no} convergence of $\sigma_N$
to $\sigma^*$.


\section{Appendix: CA symmetries}\label{symm}

There are two symmetry operations on the collection of sets of CA
rules.  The first one is an extension of the mirror map
$\mathcal{M}$ defined by
$$\mathcal{M}(0)=1,\quad \mathcal{M}(1)=0.$$
It is given by
$$(\mathcal{M}\varphi)(x,y)=\mathcal{M}(\varphi(\mathcal{M}(x),\mathcal{M}(y))).$$
The second one is space reversal $\mathcal{R}$, defined by
$$(\mathcal{R}\varphi)(x,y)=\varphi(y,x).$$
Both operations are involutions. The effect of $\mathcal{M}$, written as a permutation is:
 $$ (0\; 15) (1 \;7) (2 \;11) (4 \;13) (6\; 9) (8\; 14).$$
 The effect of $\mathcal{R}$ is:
$$ (2 \;4) (3\; 5) (10\; 12)  (11\; 13).$$

\section{ Appendix: impermeable blocks}\label{imp-ver}
Here we give two examples of the proofs that the blocks obtained
from Table~\ref{imperm} yield impermeable blocks.

To prove that the pair $(b_1,\ldots, b_4)=(0,0,1,0),
(\jfi{j_1},\ldots,\jfi{j_4})=(2,9,9,2)$ yields an impermeable
block, we have to consider a configuration as
\begin{center}
 \begin{tabular}{c|c|c|c|c|c}
\;  &  $\firule{2}$ &   $\firule{9}$ & $\firule{9}$ & $\firule{2}$ & \; \\
 \sepline
 $x$ & 0 & 0 & 1 & 0 & $y$
\end{tabular}
\end{center}
Filling in $\firule{2}(x,0)=0,\, \firule{9}(0,1)=0,\,\firule{9}(0,0)=1$
 and $\firule{2}(1,y)=0$, which is true for arbitrary $x$ and $y$,
 we obtain the next line (time $t=1$):
\begin{center}
  \begin{tabular}{c|c|c|c|c|c}
\;  &  $\firule{2}$ &   $\firule{9}$ & $\firule{9}$ & $\firule{2}$ & \; \\
 \sepline
 $x$ & 0 & 0 & 1 & 0 & $y$ \\
$x'$ & 0 & 0 & 1 & 0 & $y'$
\end{tabular}
\end{center}
Since this is the same as the state at time 0, it follows that we
do indeed have an impermeable block, consisting of four cells that
always stabilize at time 0.

There is only one exception of an impermeable block which does
\emph{not} consist of stable cells:
\begin{center}
 \begin{tabular}{c|c|c|c}
\;  &  $\firule{5}$ &   $\firule{3}$ &  \; \\
 \sepline
 $x$ & 0 & 0 & $y$\\
  $x'$ & 1 & 1 & $y'$\\
   $x''$ & 0 & 0 & $y''$\\
    $x'''$ & 1 & 1 & $y'''$\\
\end{tabular}
\end{center}
This time we obtain a period 2 impermeable block, but actually a
period 1 impermeable block is also possible by taking
$(b_1,b_2)=(0,1)$.

\section{ Appendix: no impermeable blocks}\label{no-imp}

We will show that the RBCA's with their support in
$\{\firule{2},\firule{3},\firule{11}\}$ do not admit any
impermeable blocks. We will do this by showing that diagrams as in
the previous appendix can not exist. In the following we will use
frequently that a block is impermeable if and only if its mirror
image, respectively space reversal is impermeable.


Consider the first element $\varphi_{\jfi{j_1}}$ of the
$\varphi$-block. This can not be $\firule{3}$, since
$\firule{3}(x,y)=x$, which is not compatible with impermeability
from the left. So it has to be $\firule{2}$ or $\firule{11}$.
Since $\mathcal{M}(\firule{2})=\firule{11}$, and
$\mathcal{M}(\firule{3})=\firule{3}$, we can assume without loss
of generality that it is $\firule{11}$. Since $\firule{11}(x,y)=0$
if and only if $(x,y)=(1,0)$ the first element of the iterates of
the $b$-block  must be equal to 1:
\begin{center}
 \begin{tabular}{c|c|c|c|c|c}
\;  &  $\firule{11}$ &   $\varphi_{\cdot}$ & $\varphi_{\cdot}$ & $\varphi_{\cdot}$ & \ldots \\
 \sepline
 $x$ & $\cdot$ & $\cdot$ & $\cdot$ & $\cdot$ & \ldots \\
$x'$ & $1$ & $\cdot$ & $\cdot$ & $\cdot$ & \ldots \\
$x''$ & $1$ & $\cdot$ & $\cdot$ & $\cdot$ & \ldots
\end{tabular}
\end{center}
We next consider all three possibilities for
$\varphi_{\jfi{j_2}}$, filling in, if possible,  the values of the
second cell given by the rules. For $\jfi{j_2}=2$ we obtain:
\begin{center}
 \begin{tabular}{c|c|c|c|c|c}
\; &$\firule{11}$ & $\firule{2}$&$\varphi_{\cdot}$&$\varphi_{\cdot}$&\ldots\\
\sepline
 $x$ & $\cdot$ & $\cdot$ & $\cdot$ & $\cdot$ & $\ldots$ \\
 $x'$ & $1$ & $\cdot$ & $\cdot$ & $\cdot$ &  $\ldots$ \\
 $x''$ & $1$ & $0$ & $\cdot$ & $\cdot$ & $\ldots$
\end{tabular}
\end{center}
However this gives a contradiction at $t=3$, since
$\firule{11}(x'',0)=1-x''$ depends on $x''$. Exactly the same
contradiction occurs when $\jfi{j_2}=3$. Consequently we must have
$\jfi{j_2}=11$. The whole second column of the $b$-block must
consist of 1's, again because of $\firule{11}(x,0)=1-x$. But then,
since $\firule{11}(1,0)=0$ also the 3$^{rd}$ column must be filled
with 1's starting from $t=1$. It is then quickly checked in the
following diagram that $\jfi{j_3}=2$ and $\jfi{j_3}=3$ are
impossible:
\begin{center}
 \begin{tabular}{c|c|c|c|c|c}
\; &$\firule{11}$ & $\firule{11}$&$\varphi_{\jfi{j_3}}$&$\varphi_{\cdot}$&\ldots\\
\sepline
 $x$ & $\cdot$ & $1$ & $\cdot$ & $\cdot$ & $\ldots$ \\
 $x'$ & $1$ & $1$ & $1$ & $\cdot$ &  $\ldots$ \\
 $x''$ & $1$ & $1$ & $1$ & $\cdot$ & $\ldots$
\end{tabular}
\end{center}
Conclusion: also $\jfi{j_3}=11$. Continuing in this fashion we
find that the (string of indices of the) $\varphi$-block has the
form $(\jfi{j_1},\ldots,\jfi{j_p})=(11,\ldots,11)$. But then we
have a problem in the $p^{\rm th}$ column since
$\firule{11}(1,y)=y$ depends on $y$. Hence for any $p$ this
possibility is ruled out.

\nocite{*}

\bibliographystyle{plain}

\bibliography{RBCA}


\end{document}